\documentclass[12pt]{article}                                               

\input{amssym.def}                                                           
\input{amssym.tex}

\begin{document}                                                             
\title{Cluster $C^*$-algebras and knot polynomials}

\author{Igor  ~Nikolaev
}


\date{}
 \maketitle


\newtheorem{thm}{Theorem}
\newtheorem{lem}{Lemma}
\newtheorem{dfn}{Definition}
\newtheorem{rmk}{Remark}
\newtheorem{cor}{Corollary}
\newtheorem{cnj}{Conjecture}
\newtheorem{exm}{Example}


\newcommand{\Coh}{\hbox{\bf Coh}}
\newcommand{\Mod}{\hbox{\bf Mod}}
\newcommand{\Tors}{\hbox{\bf Tors}}

\begin{abstract}
We construct  a representation of the braid groups in a cluster $C^*$-algebra
coming from a triangulation of the Riemann surface  $S$ with 
one or two cusps.  It is shown that the Laurent 
polynomials  attached to the $K$-theory of such an algebra 
are  topological invariants  of  the closure of braids.  In particular, 
 the Jones and HOMFLY polynomials of a  knot  
correspond to the case $S$ being a sphere  with two cusps and a torus  
with one cusp, respectively.

\vspace{7mm}


\vspace{5mm}

\end{abstract}

\section{Introduction}
The {\it trace invariant} $V_L(t)$ of  a link $L$  was introduced in [Jones 1985]  \cite{Jon1}.
The $V_L(t)$ is a Laurent polynomial $\mathbf{Z}[t^{\pm {1\over 2}}]$ obtained from 
a representation  of the braid group $B_k$  in  an operator (von Neumann) algebra $A_k$;  
 the canonical trace  on $A_k$  times a multiple  ${t\over (1+t)^2}$  is invariant of the Markov 
move of type II of a braid $b\in B_k$.    The Jones polynomial  $V_L(t)$ is a powerful 
topological invariant of a  link  $L$  obtained by the closure of  $b$. The algebra $A_k$ itself comes from 
an analog of  the Galois theory for   von Neumann algebras  called  {\it basic construction}
 [Jones 1991, Section 2.6]   \cite{J1}.
It is yet unclear why the braid groups  appear   in  the context  of  operator algebras, 
 let alone  an invariant  trace with such a remarkable property;  likewise, 
one can seek to extend the invariant  $V_L(t)$
to the multivariable Laurent  polynomials.

{\it Cluster algebras}  are a class of commutative rings introduced by 
[Fomin \& Zelevinsky 2002]  \cite{FoZe1}  having deep roots  in    
hyperbolic  geometry and Teichm\"uller theory  [Williams 2014] \cite{Wil1}.  
Such an algebra  ${\cal A}(\mathbf{x}, B)$ 
is a subring of the field  of  rational functions in $n$ variables
depending  on a cluster  $\mathbf{x}=(x_1,\dots, x_n;  ~y_1, \dots, y_m)$
of {\it mutable} variables $x_i$ and {\it frozen} variables $y_i$
and a skew-symmetric matrix  $B=(b_{ij})\in M_n(\mathbf{Z})$; 
the pair  $(\mathbf{x}, B)$ is called a {\it seed}.
In terms of the coefficients $c_i$  from  a semi-field 
$({\Bbb P}, \oplus, \bullet)$ a new cluster 
$\mathbf{x}'=(x_1,\dots,x_k',\dots,  x_n;
~c_1,\dots,c_j',\dots, c_n)$ and a new
skew-symmetric matrix $B'=(b_{ij}')$ is obtained from 
$(\mathbf{x}, B)$ by the  exchange relations:
\begin{eqnarray}\label{eq1}
b_{ij}' &=& \cases{-b_{ij}  & \hbox{if}  $i=k$ \hbox{or} $j=k$\cr
b_{ij}+{|b_{ik}|b_{kj}+b_{ik}|b_{kj}|\over 2}  & \hbox{otherwise,}}\cr
c_j' &=& \cases{{1\over c_k}  & \hbox{if}  $j=k$\cr
{c_jc_k^{\max(b_{kj}, 0)}\over  (c_k\oplus 1)^{b_{kj}}} & \hbox{otherwise,}}\cr
x_k'  &=&{c_k \prod_{i=1}^n  x_i^{\max(b_{ik}, 0)} + \prod_{i=1}^n  x_i^{\max(-b_{ik}, 0)}\over 
(c_k\oplus 1) ~x_k}, 
\end{eqnarray}
see  [Williams 2014, Definition 2.22]  \cite{Wil1} for the details.  
The seed $(\mathbf{x}', B')$ is said to be a {\it mutation} of $(\mathbf{x}, B)$ in direction $k$,
where $1\le k\le n$;   the algebra  ${\cal A}(\mathbf{x}, B)$ is  generated by cluster  variables $\{x_i\}_{i=1}^{\infty}$
obtained from the initial seed $(\mathbf{x}, B)$ by the iteration of mutations  in all possible
directions $k$.   The {\it Laurent phenomenon} 
proved by  [Fomin \& Zelevinsky 2002]  \cite{FoZe1}  says  that  ${\cal A}(\mathbf{x}, B)\subset \mathbf{Z}[\mathbf{x}^{\pm 1}]$,
where  $\mathbf{Z}[\mathbf{x}^{\pm 1}]$ is the ring of  the Laurent polynomials in  variables $\mathbf{x}=(x_1,\dots,x_n)$;  
in other words, each  generator $x_i$  of  algebra ${\cal A}(\mathbf{x}, B)$  can be 
written as a  Laurent polynomial in $n$ variables with the   integer coefficients. 
(Note that  the Laurent phenomenon turns ${\cal A}(\mathbf{x}, B)$ into an additive abelian 
group with an order coming from the semi-group of the Laurent polynomials with positive
coefficients.)     
In what follows,   we deal with a cluster algebra  ${\cal A}(\mathbf{x},  S_{g,n})$
coming from  a triangulation of the Riemann surface  $S_{g,n}$   of genus $g$  with $n$ cusps,
see  [Fomin,  Shapiro  \& Thurston  2008]  \cite{FoShaThu1} the details.

{\it Cluster $C^*$-algebras}  are a class of non-commutative rings  ${\Bbb A}(\mathbf{x}, B)$,
such that $K_0({\Bbb A}(\mathbf{x}, B))\cong {\cal A}(\mathbf{x}, B)$  \cite{Nik1};  
here $K_0({\Bbb A}(\mathbf{x}, B))$ is the $K_0$-group of a $C^*$-algebra ${\Bbb A}(\mathbf{x}, B)$
and $\cong$ is an isomorphism of the additive abelian groups with order [Blackadar 1986] \cite{B}.  
The ${\Bbb A}(\mathbf{x}, B)$  is an {\it $AF$-algebra}
 given by the Bratteli diagram [Bratteli 1972]  \cite{Bra1}; 
 such a diagram can be obtained from a  mutation tree of the initial seed $(\mathbf{x}, B)$
 modulo an equivalence relation between the seeds  lying at  the same level.
It is known that the $AF$-algebras are   characterized 
by their $K$-theory [Elliott 1976]   \cite{Ell1}.  
Equivalently,  the ${\Bbb A}(\mathbf{x}, B)$ is an algebra
over the complex numbers  generated by a series of projections 
$\{e_i\}_{i=1}^{\infty}$.

The aim of our note is  a representation of the braid group 
 $B_k=$\linebreak
 $\{\sigma_1,\dots,\sigma_{k-1} ~|~ \sigma_i\sigma_{i+1}\sigma_i=
 \sigma_{i+1}\sigma_i\sigma_{i+1},  ~\sigma_i\sigma_j=\sigma_j\sigma_i
 ~\hbox{if}  ~|i-j|\ge 2\}$
 into   an algebra ${\Bbb A}(\mathbf{x}, S_{g,n})$,  so that  the Laurent  phenomenon in 
$K_0({\Bbb A}(\mathbf{x}, S_{g,n}))$  corresponds  to the polynomial  invariants of
 the closure  of  braids   $b\in B_k$.  
In particular, if $g=0$ and $n=2$ or $g=n=1$ one recovers 
the  Jones invariant  $V_L(t)$ or  the HOMFLY  polynomials   of knots [Freyd, Yetter, Hoste, Lickorish, Millet \&
 Ocneanu 1985] \cite{FYHLMO},  respectively.    Whenever  $g\ge 2$  one gets 
 new topological invariants   generalizing   the Jones and HOMFLY
polynomials to an arbitrary (but finite) number of variables. 
  The  $AF$-algebra ${\Bbb A}(\mathbf{x}, S_{g,n})$  itself  can be viewed 
as an analog of the tower $\cup_{k=1}^{\infty}  A_k$  of von Neumann algebras 
$A_k$ arising from   the basic construction   [Jones 1991, Section 3.4]   \cite{J1}. 
Unlike [Jones 1985]  \cite{Jon1},  we  exploit the  phenomenon of cluster algebras 
and the {\it Birman-Hilden  Theorem}  relating the  braid groups $B_{2g+n}$ with the mapping class group 
 of surface $\{S_{g,n} ~|~ n=1; 2\}$  [Birman \& Hilden 1971]  \cite{BiHi1}.
 Our main result can be formulated   as follows. 
\begin{thm}\label{thm1}
The formula $\sigma_i\mapsto e_i+1$ defines a representation 
\begin{equation}
\rho: \cases{B_{2g+1}\to {\Bbb A}(\mathbf{x}, S_{g,1}) &\cr
B_{2g+2}\to {\Bbb A}(\mathbf{x}, S_{g,2}). &}
\end{equation}
If  $b\in B_{2g+1}$  ($b\in B_{2g+2}$, resp.) is a braid,   there exists 
a Laurent polynomial  $[\rho(b)]\in   K_0({\Bbb A}(\mathbf{x}, S_{g,1}))$   
($[\rho(b)] \in K_0({\Bbb A}(\mathbf{x}, S_{g,2}))$, resp.)
with the integer coefficients    depending on  $2g$  ($2g+1$, resp.) 
variables,  such that  $[\rho(b)]$  is a topological invariant of the closure of 
  $b$.  
\end{thm}
The article is organized as follows.  We introduce preliminary facts and notation 
in Section 2.   Theorem \ref{thm1} is proved in Section 3. 
To illustrate theorem \ref{thm1} in Section 4,   we consider  the cases  $g=0, ~n=2$ and $g=n=1$
corresponding  to the Jones and HOMFLY polynomials, respectively.

\section{Preliminaries}
We shall briefly review the Birman-Hilden Theorem,  the cluster $C^*$-algebras
and knots. We refer the reader to   [Birman \& Hilden 1971]  \cite{BiHi1},
[Blackadar 1986] \cite{B}, 
  [Farb \& Margalit  2011]   \cite{FM},  
 [Fomin,  Shapiro  \& Thurston  2008]  \cite{FoShaThu1},    [Jones 1985]  \cite{Jon1},
[Jones 1991,  Lecture 6]   \cite{J1},   [Williams 2014] \cite{Wil1},  \cite{Nik2} and \cite{Nik1}  for a detailed account.

\subsection{Birman-Hilden Theorem}
Let $S_{g,n}$ be a Riemann surface  of genus $g\ge 0$ with $n\ge 1$  cusps and  such that $2g-2+n>0$;
denote by   $T_{g,n}\cong {\Bbb R}^{6g-6+2n}$   the (decorated)  Teichm\"uller space of $S_{g,n}$,
i.e.  a collection   of all Riemann surfaces of genus  $g$ with $n$ cusps endowed with the natural topology.
By $Mod~S_{g,n}$ we understand the mapping class group of surface $S_{g,n}$, i.e.
a group of the homotopy classes of  orientation preserving automorphisms of $S_{g,n}$ fixing all cusps. 
It is well known that two Riemann surfaces $S,S'\in T_{g,n}$ are isomorphic 
if and only if there exists a $\varphi\in   Mod~S_{g,n}$,  such that
$S'=\varphi(S)$; thus each $\varphi\in   Mod~S_{g,n}$ corresponds to a homeomorphism
of the Teichm\"uller space $T_{g,n}$.  
If $\gamma$ is a simple closed curve on $S_{g,n}$,  let $D_{\gamma}\in   Mod~S_{g,n}$ 
be the Dehn twist around $\gamma$;  a pair of the Dehn twists ${D_{\gamma_i}}$ and  ${D_{\gamma_j}}$ 
satisfy the braid relations:
\begin{equation}\label{eq3}
\cases{D_{\gamma_i} D_{\gamma_j}  D_{\gamma_i}=D_{\gamma_j} D_{\gamma_i} D_{\gamma_j},  
& if   ~$\gamma_i\cap\gamma_j=\{\hbox{\sf single point}\}$ \cr
D_{\gamma_i}D_{\gamma_j}=D_{\gamma_j} D_{\gamma_i},  & if  ~$\gamma_i\cap\gamma_j=\emptyset.$
}
\end{equation}
A system of simple closed curves $\{\gamma_i\}$ on $S_{g,n}$ is called a {\it chain},
if 
\begin{equation}
\cases{\gamma_i\cap\gamma_{i+1}=\{\hbox{\sf single point}\}
&\cr
\gamma_i\cap\gamma_j=\emptyset   \quad\hbox{otherwise.}
&
}
\end{equation}
Consider a chain $\{\gamma_1,\dots,\gamma_{2g+1}\}$ shown in 
[Farb \& Margalit  2011,  Figure 2.7]   \cite{FM}.  
The fundamental domain of $S_{g,n}$ obtained by a cut along the chain
is a $(4g+2)$-gon with the opposite sides identified [Farb \& Margalit  2011,  Figure 2.2]   \cite{FM}.
A {\it hyperelliptic involution} $\iota\in   Mod~S_{g,n}$ is a rotation by the angle $\pi$ of the 
$(4g+2)$-gon;  clearly, the chain    $\{\gamma_1,\dots,\gamma_{2g+1}\}$ is an invariant
of the involution $\iota$.

In what follows, we focus on the case $n=1$ or $2$.    Notice that $S_{g,1}$ ($S_{g,2}$, resp.) 
can be replaced by a surface of the same genus with one (two, resp.) boundary components
and no cusps;  since the mapping class group preserves the boundary components,  we 
work with this new surface while keeping the old notation. 
It is known that $S_{g,1}$ ($S_{g,2}$, resp.)  is a double cover of a disk ${\cal D}_{2g+1}$
(${\cal D}_{2g+2}$, resp.)  ramified at the $2g+1$ ($2g+2$, resp.) inner points of the disk;  
we refer  the reader to    [Farb \& Margalit  2011,  Figure 9.15]   \cite{FM} for a picture. 
It transpires that each automorphism of ${\cal D}_{2g+1}$  (${\cal D}_{2g+2}$, resp.)
pulls back to an automorphism of  $S_{g,1}$ ($S_{g,2}$, resp.) commuting with the 
hyperelliptic involution $\iota$.  Recall that  $Mod ~{\cal D}_{2g+1}\cong B_{2g+1}$  
($Mod ~{\cal D}_{2g+2}\cong B_{2g+2}$, resp.);  a subgroup of $Mod ~S_{g,1}$ ($Mod ~S_{g,2}$, resp.)
commuting with $\iota$ is called {\it symmetric} and denoted by  $SMod ~S_{g,1}$ ($SMod ~S_{g,2}$, resp.)
The following result is critical. 
\begin{thm}\label{thm2}
{\bf ([Birman \& Hilden 1971]  \cite{BiHi1})}
There exists an isomorphism: 
\begin{equation}
\cases{B_{2g+1}\cong  SMod ~S_{g,1} &\cr
B_{2g+2}\cong  SMod ~S_{g,2} &}
\end{equation}
given by the formula $\sigma_i\mapsto D_{\gamma_i}$,  where $\sigma_i$
is a generator of the braid group $B_{2g+1}$ ($B_{2g+2}$, resp.)  
and $D_{\gamma_i}$ is the Dehn twist around the simple closed curve 
$\gamma_i$ of a  chain in the  $S_{g,1}$ ($S_{g,2}$, resp.)
\end{thm}

\subsection{Cluster $C^*$-algebras from Riemann surfaces}
The fundamental domain of the Riemann surface  $S_{g,n}$ has a  triangulation by
the  $6g-6+3n$ geodesic arcs $\gamma_i$;   the endpoint  of  each $\gamma_i$
is a cusp at  the absolute of Lobachevsky plane ${\Bbb H}=\{x+iy\in {\Bbb C} ~|~ y>0\}$.
Denote by  $l(\gamma_i)$  the $\pm$ 
hyperbolic length  of $\gamma_i$  between two horocycles
around the endpoints of $\gamma_i$;   consider  the $\lambda(\gamma_i)=\exp ~({1\over 2} l(\gamma_i))$. 
The following result says that  $\lambda(\gamma_i)$ are coordinates in the Teichm\"uller space 
$T_{g,n}$.  
\begin{thm}\label{thm3}
{\bf ([Penner 1987]  \cite{Pen1})}  
The map $\lambda: \{\gamma_i\}_{i=1}^{6g-6+3n}\to T_{g,n}$ 
is  a homeomorphism.
 \end{thm}
\begin{rmk}
\textnormal{
The six-tuples of numbers  $\lambda(\gamma_i)$  must satisfy the {\it Ptolemy relation}
$\lambda(\gamma_1)\lambda(\gamma_2)+\lambda(\gamma_3)\lambda(\gamma_4)=
\lambda(\gamma_5)\lambda(\gamma_6)$,
where $\gamma_1, \dots, \gamma_4$ are pairwise opposite sides and $\gamma_5, \gamma_6$
are the diagonals of a geodesic quadrilateral in ${\Bbb H}$.
The $n$ Ptolemy relations reduce the number of independent variables $\lambda(\gamma_i)$ to 
$6g-6+2n=\dim ~T_{g,n}$.
}
\end{rmk}
Let $T=\{\gamma_i\}_{i=1}^{6g-6+3n}$ be a triangulation of  $S_{g,n}$; 
consider a skew-symmetric matrix $B_T=(b_{ij})$ of rank $6g-6+3n$,  where  $b_{ij}$ is equal to the number 
of triangles in $T$ with sides $\gamma_i$ and $\gamma_j$ in clockwise order minus 
the number of triangles in $T$ with  sides $\gamma_i$ and $\gamma_j$ in the 
counter-clockwise order. (For a quick example of matrix $B_T$ we refer the reader to
Section 4.) 
\begin{thm}\label{thm4}
{\bf ([Fomin, Shapiro \& D.~Thurston 2008]  \cite{FoShaThu1})}  
The cluster algebra ${\cal A}(\mathbf{x}, B_T)$ does not depend on triangulation
$T$, but only on the surface $S_{g,n}$;  namely,   replacement of the geodesic 
arc $\gamma_k$ by a new geodesic arc $\gamma_k'$ (a flip of $\gamma_k$) 
corresponds to a mutation $\mu_k$ of the seed $(\mathbf{x}, B_T)$. 
\end{thm}
\begin{dfn}
The algebra ${\cal A}(\mathbf{x}, S_{g,n}):={\cal A}(\mathbf{x}, B_T)$ 
is called a cluster algebra of the Riemann surface $S_{g,n}$.   
\end{dfn}
A {\it $C^*$-algebra} is an algebra $A$ over $\mathbf{C}$ with a norm
$a\mapsto ||a||$ and an involution $a\mapsto a^*$ such that
it is complete with respect to the norm and $||ab||\le ||a||~ ||b||$
and $||a^*a||=||a^2||$ for all $a,b\in A$.
Any commutative $C^*$-algebra is  isomorphic
to the algebra $C_0(X)$ of continuous complex-valued
functions on some locally compact Hausdorff space $X$; 
otherwise, $A$ represents a noncommutative  topological
space.   For a unital $C^*$-algebra $A$, let $V(A)$
be the union (over $n$) of projections in the $n\times n$
matrix $C^*$-algebra with entries in $A$;
projections $p,q\in V(A)$ are (Murray - von Neumann)  {\it equivalent} if there exists a partial
isometry $u$ such that $p=u^*u$ and $q=uu^*$. The equivalence
class of projection $p$ is denoted by $[p]$;
the equivalence classes of orthogonal projections can be made to
a semigroup by putting $[p]+[q]=[p+q]$. The Grothendieck
completion of this semigroup to an abelian group is called
the  $K_0$-group of the algebra $A$.
The functor $A\to K_0(A)$ maps the category of unital
$C^*$-algebras into the category of abelian groups, so that
projections in the algebra $A$ correspond to a positive
cone  $K_0^+\subset K_0(A)$ and the unit element $1\in A$
corresponds to an order unit $u\in K_0(A)$.
The ordered abelian group $(K_0,K_0^+,u)$ with an order
unit  is called a {\it dimension group};  an order-isomorphism
class of the latter we denote by $(G,G^+)$.

An {\it $AF$-algebra}  ${\Bbb A}$  (Approximately Finite $C^*$-algebra) is defined to
be the  norm closure of an ascending sequence of   finite dimensional
$C^*$-algebras $M_n$,  where  $M_n$ is the $C^*$-algebra of the $n\times n$ matrices
with entries in $\mathbf{C}$;  such an algebra is given by an infinite graph called 
{\it Bratteli diagram},  see  [Bratteli 1972]  \cite{Bra1} for a definition. 
The dimension group
$(K_0({\Bbb A}), K_0^+({\Bbb A}), u)$ is a complete isomorphism
invariant of the algebra ${\Bbb A}$ [Elliott 1976]   \cite{Ell1}.  
The order-isomorphism  class $(K_0({\Bbb A}), K_0^+({\Bbb A}))$
 is an invariant of the {\it Morita equivalence} of algebra 
${\Bbb A}$,  i.e.  an isomorphism class in the category of 
finitely generated projective modules over ${\Bbb A}$.    
The {\it scale} $\Gamma$ is a subset of  $K_0^+({\Bbb A})$
which is generating, hereditary and directed, i.e. 
(i)  for each $a\in K_0^+({\Bbb A})$ there exist $a_1,\dots,a_r\in\Gamma({\Bbb A})$,
such that $a=a_1+\dots+a_r$;
(ii)  if $0\le a\le b\in \Gamma$, then $a\in\Gamma$;
(iii) given $a,b\in\Gamma$ there exists $c\in\Gamma$, such that
$a,b\le c$.  If $u$ is an order unit,  then the set $\Gamma:=\{a\in K_0^+({\Bbb A})~|~ 0\le a\le u\}$
is a scale;  thus the dimension group of algebra ${\Bbb A}$ can be written 
in the form $(K_0({\Bbb A}), K_0^+({\Bbb A}), \Gamma)$. 
\begin{dfn}
{\bf (\cite{Nik1})}
By a cluster $C^*$-algebra ${\Bbb  A}(\mathbf{x}, S_{g,n})$ one understands an $AF$-algebra 
satisfying an isomorphism of the scaled  dimension groups:
\begin{equation}
(K_0({\Bbb  A}(\mathbf{x}, S_{g,n})), K_0^+({\Bbb  A}(\mathbf{x}, S_{g,n})), u) \cong
({\cal   A}(\mathbf{x}, S_{g,n}), {\cal   A}^+(\mathbf{x}, S_{g,n}), u'), 
\end{equation}
where ${\cal   A}^+(\mathbf{x}, S_{g,n})$  is a semi-group of the Laurent polynomials with positive 
coefficients and $u'$ is an order unit in ${\cal   A}^+(\mathbf{x}, S_{g,n})$.  
\end{dfn}
\begin{rmk}\label{rmk2}
\textnormal{
Theorems \ref{thm3} and \ref{thm4} imply that   the algebra  ${\Bbb  A}(\mathbf{x}, S_{g,n})$
is a non-commutative coordinate ring of the Teichm\"uller space  $T_{g,n}$;  in other words,  the 
diagram in Figure 1 must be commutative.  
}
\end{rmk}
\begin{figure}[here]
\begin{picture}(300,110)(-100,-5)
\put(20,70){\vector(0,-1){35}}
\put(130,70){\vector(0,-1){35}}
\put(52,23){\vector(1,0){53}}
\put(52,83){\vector(1,0){53}}
\put(-5,20){${\Bbb  A}(\mathbf{x}, S_{g,n})$}
\put(115,20){${\Bbb  A}(\mathbf{x}, S_{g,n})$}
\put(10,80){$T_{g,n}$}
\put(120,80){$T_{g,n}$}
\put(30, 95){homeomorphism}
\put(65, 30){inner}
\put(45, 5){automorphism}
\end{picture}
\caption{Coordinate ring of the space  $T_{g,n}$.}
\end{figure}
\begin{rmk}\label{rmk3}
\textnormal{
The ${\Bbb  A}(\mathbf{x}, S_{g,n})$ is the  norm-closure of  an  algebra 
of the non-commutative polynomials $\mathbf{C}\langle e_1,e_2,\dots\rangle$,
 where   $\{e_i\}_{i=1}^{\infty}$ are projections  in the algebra  ${\Bbb  A}(\mathbf{x}, S_{g,n})$;  
this fact follows from the $K$-theory of   ${\Bbb  A}(\mathbf{x}, S_{g,n})$.  
On the other hand, the algebra ${\cal  A}(\mathbf{x}, S_{g,n})$ is generated 
by the cluster variables $\{x_i\}_{i=1}^{\infty}$.  We shall denote by $\rho(x_i)=e_i$ 
a natural bijection between the two sets of generators. 
}
\end{rmk}

\subsection{Knots and links}
A {\it knot}  is a tame embedding of the circle $S^1$  into the Euclidean 
space $\mathbf{R}^3$;  a {\it link} with the $n$ components is such an
embedding of the union $S^1\cup\dots\cup S^1$.  The classification 
of (distinct)  knots and links is a difficult open problem of topology. 
The {\it Alexander Theorem} says that every knot or link comes from 
the closure of a braid $b\in B_k=
\{\sigma_1,\dots,\sigma_{k-1} ~|~ \sigma_i\sigma_{i+1}\sigma_i=
 \sigma_{i+1}\sigma_i\sigma_{i+1},  ~\sigma_i\sigma_j=\sigma_j\sigma_i
 ~\hbox{if}  ~|i-j|\ge 2\}$,
 i.e. tying the top end of each string of $b$ to the end of a string in the 
 same position at the bottom.  Thus the braids can classify knots and 
 links but sadly rather ``unrelated'' braids $b\in B_k$ and $b'\in B_{k'}$
 can produce the same knot.  Namely, the {\it Markov Theorem}  says 
 that the closure of braids  $b,b'\in B_k$ corresponds to the same knot 
 or link,  if and only if: (i) $b'=gbg^{-1}$ for a braid $g\in B_k$ or
 (ii)  $b'=b\sigma_k^{\pm 1}$ for the generator $\sigma_k\in B_{k+1}$.
 Thus the {\it Markov move of type II}   always pushes ``to infinity'' 
 the desired classification hinting that an asymptotic invariant 
 is required;  below we consider two examples of such invariants.

 \medskip
 The {\it Jones polynomial} of the closure $L$ of a braid $b\in B_k$ 
 is defined by the formula:
\begin{equation}
V_L(t)=\left(-{t+1\over\sqrt{t}}\right)^{k-1}~tr~(r_t(b)),
\end{equation} 
where  $r_t$ is a representation of $B_k$ in a 
von-Neumann algebra $A_k$ and $tr$ is a  trace function;   
 the $V_L(t)\in\mathbf{Z}[t^{\pm {1\over 2}}]$, i.e. a Laurent
 polynomial in the variable  $\sqrt{t}$ [Jones 1985]  \cite{Jon1}. 
 If $K$ is the unknot then $V_K(t)=1$ and each polynomial $V_L(t)$ can 
 be calculated  from $K$ using the {\it skein relation}:
\begin{equation}\label{skein}
{1\over t}V_{L^-} -tV_{L^+}=\left(\sqrt{t}-{1\over\sqrt{t}}\right)V_L, 
\end{equation} 
where $L^+$ ($L^-$, resp.) is a link obtained by adding an overpass
(underpass, resp.) to the link $L$ [Jones 1985, Theorem 12]  \cite{Jon1}.

\medskip 
The {\it HOMFLY polynomial}  of a link $L$ is a Laurent polynomial 
$\rho_L(l,m)\in\mathbf{Z}[l^{\pm 1}, m^{\pm 1}]$;   the  $\rho_L(l,m)$
is defined recursively from the HOMFLY polynomial 
$\rho_K(l,m)=1$ of  the  unknot $K$ using the skein relation:  
\begin{equation}\label{homfly}
l\rho_{L^+} +{1\over l}\rho_{L^-}+ m\rho_L=0, 
\end{equation} 
where $L^+$ ($L^-$, resp.) is a link obtained by adding an overpass
(underpass, resp.) to   $L$  [Freyd, Yetter, Hoste, Lickorish, Millet \&
 Ocneanu 1985,  Remark 3]  \cite{FYHLMO}.

\section{Proof of theorem \ref{thm1}}
For the sake of clarity, let us  outline  the main ideas.  Roughly speaking, 
the Birman-Hilden's  Theorem \ref{thm2}   says that a generator $\sigma_i\in B_{2g+1}$ 
($B_{2g+2}$, resp.)  is given by  the Dehn twist $D_{\gamma_i}\in
Mod~S_{g,1}$  ($Mod~S_{g,2}$, resp.) around a closed curve $\gamma_i$.  
The $D_{\gamma_i}$ itself  is a homeomorphism  of the Teichm\"uller space  
$T_{g,1}$  ($T_{g,2}$, resp.);  therefore the $D_{\gamma_i}$  induced  an inner 
automorphism  $x\mapsto u_i xu^{-1}_i$, where $u_i$ is a unit  of  the cluster 
$C^*$-algebra ${\Bbb  A}(\mathbf{x}, S_{g,1})$ (${\Bbb  A}(\mathbf{x}, S_{g,2})$, resp.)
(We refer the reader to  Figure 1.)  
On the other hand,  the Fomin-Shapiro-D.~Thurston's  Theorem \ref{thm4} and remark \ref{rmk3}
imply that units  $u_i$ and projections  $e_i$ in  algebra ${\Bbb  A}(\mathbf{x}, S_{g,1})$  
(${\Bbb  A}(\mathbf{x}, S_{g,2})$, resp.) are bijective. But  the minimal degree polynomial  
in variable  $e_i$  corresponding to  a unit $u_i$ is a linear polynomial of  the form $u_i=ae_i+b$, 
where  $a$ and $b$ are complex constants.   (The inverse is given by the formula
$u_i^{-1}=-{a\over (a+b)b}e_i +{1\over b}$.)  
The Birman-Hilden's  Theorem \ref{thm2}  implies that  the  $u_i$ satisfy the braid relations 
$\{u_i u_{i+1} u_i= u_{i+1} u_i u_{i+1},  ~u_i u_j=u_j u_i,  ~\hbox{if}  ~|i-j|\ge 2\}$;
moreover,  if one substitutes  $u_i=ae_i+b$ in the braid relations,
then:  
\begin{equation}\label{eq8}
\cases{ e_i^2=e_i,  &\cr
e_ie_{i\pm 1}e_i=-{(a+b)b\over a^2}e_i, &\cr
e_ie_j=e_je_i,  & if  $|i-j|\ge 2.$}
\end{equation}
\begin{rmk}\label{TL}
\textnormal{
The  relations (\ref{eq8}) are invariant of the involution $e_i^*=e_i$ if and only if  
 ${(a+b)b\over a^2}\in\mathbf{R}$;  in this case  the algebra  ${\Bbb  A}(\mathbf{x}, S_{g,1})$  
 (${\Bbb  A}(\mathbf{x}, S_{g,2})$, resp.)
contains  a finite-dimensional $C^*$-algebra ${\Bbb A}_{2g}$ (${\Bbb A}_{2g+1}$, resp.) 
obtained from  the norm closure of a self-adjoint representation  of a  
 {\it Temperley-Lieb algebra}.
 }
 \end{rmk}    
 Thus one gets  a representation $\rho:  B_{2g+1}\to {\Bbb  A}(\mathbf{x}, S_{g,1})$ 
 ($B_{2g+2}\to {\Bbb  A}(\mathbf{x}, S_{g,2})$, resp.)
given by the formula $\sigma_i\mapsto ae_i+b$,  where $1\le i\le 2g$ ($1\le i\le 2g+1$, resp.)
It follows from relations (\ref{eq8})  that  the set ${\cal E} :=$\linebreak
$\{(e_{i_1}e_{i_1-1}\dots e_{j_1}) \dots  (e_{i_p}e_{i_p-1}\dots e_{j_p}) ~|~ 
1  \le  i_1<\dots <i_p<2g  ~(2g+1, \hbox{resp.}); $\linebreak
$~1 \le  j_1<\dots <j_p<2g ~(2g+1, \hbox{resp.}); 
~j_1\le i_1, \dots, j_p\le i_p\}$  is multiplicatively closed; 
moreover,  $|{\cal E}|\le {1\over n+1}\left(\small\matrix{2n\cr n}\right)=n$'th  Catalan  number,
where $n=2g$ ($n=2g+1$, resp.)
In particular,  the ${\Bbb A}_{2g}$ (${\Bbb A}_{2g+1}$, resp.)  of remark \ref{TL} is a finite-dimensional
$C^*$-algebra 
and each element $\varepsilon\in {\cal E}$ is equivalent to a projection;  the (Murray-von Neumann) equivalence 
class of the projection  will be denoted by $[\varepsilon]$.  
If $\{b=\sigma_1^{k_1}\dots\sigma_{n-1}^{k_{n}}\in B_{n+1} ~|~ k_i\in\mathbf{Z}\}$  is a braid for  
$n=2g$ ($n=2g+1$, resp.)  then the polynomial $\rho(b)=(ae_1+b)^{k_1}\dots(ae_{n-1}+b)^{k_{n}}$
unfolds into a finite sum 
$\left\{\sum_{i=1}^{|{\cal E}|} a_i\varepsilon_i  ~|~ \varepsilon_i\in {\cal E}, ~a_i\in \mathbf{Z}\right\}$. 
Therefore one gets an inclusion  
$[\rho(b)] := \left\{\sum_{i=1}^{|{\cal E}|} a_i [\varepsilon_i]  ~|~ [\varepsilon_i]\in K_0({\Bbb A}(\mathbf{x}, S_{g,1})), 
~a_i\in \mathbf{Z}\right\}\in  K_0 ({\Bbb A}(\mathbf{x}, S_{g,1})$ ($K_0 ({\Bbb A}(\mathbf{x}, S_{g,2})$, resp.) 
But  $K_0({\Bbb A}(\mathbf{x}, S_{g,n}))\cong {\cal  A}(\mathbf{x}, S_{g,n})\subset 
\mathbf{Z}[\mathbf{x}^{\pm 1}]$,  where $n=1$ ($n=2$, resp.);  thus  $[\rho(b)]\in K_0 ({\Bbb A}(\mathbf{x}, S_{g,n})$
is a Laurent polynomial with the integer coefficients  depending on  $2g$ ($2g+1$, resp.) 
variables.  
The $[\rho(b)]$  is in fact  a topological invariant of the closure of  $b$.
Indeed, $[\rho(gbg^{-1})]=[\rho(b)]$ for all $b\in B_{2g+1}$ ($B_{2g+2}$, resp.) 
because the $K_0$-group and a canonical trace $\tau$ on the algebra  ${\Bbb A}(\mathbf{x}, S_{g,1})$ 
(${\Bbb A}(\mathbf{x}, S_{g,2})$, resp.) are related [Blackadar 1986,  Section 7.3]  \cite{B};  
since the $\tau$ is a character of the representation $\rho$,  one gets the formula  $[\rho(gbg^{-1})]=[\rho(b)]$.  
An invariance of the $[\rho(b)]$  with respect to  the Markov move of type II  is a bit  subtler,
but  follows from a stability  of the $K$-theory   [Blackadar 1986, Section 5.1]  \cite{B}.
Namely, the map $\sigma_{2g+1}^{\pm 1}\mapsto 2e_{2g+1}-1$   
($\sigma_{2g+2}^{\pm 1}\mapsto 2e_{2g+2}-1$, resp.)  gives rise to a
crossed product $C^*$-algebra ${\Bbb A}_{2g}\rtimes_{\alpha} G$ (${\Bbb A}_{2g+1}\rtimes_{\alpha} G$, resp.),
where $G\cong \mathbf{Z}/2\mathbf{Z}$;
the crossed product is isomorphic to the algebra $M_2({\Bbb A}_{2g}^{\alpha})$ 
($M_2({\Bbb A}_{2g+1}^{\alpha})$, resp.), where ${\Bbb A}_{2g}^{\alpha}$
(${\Bbb A}_{2g+1}^{\alpha}$, resp.)   is the fixed-point algebra of the automorphism $\alpha$
[Fillmore 1996, Section 3.8.5]  \cite{F}.  But $K_0(M_2({\Bbb A}_{2g}^{\alpha}))\cong K_0({\Bbb A}_{2g}^{\alpha})$
by the stability of the $K$-theory;  the crossed product itself consists of the  formal sums     
$\sum_{\gamma\in G} a_{\gamma} u_{\gamma}$ and one easily derives that 
$[\rho(b\sigma_{2g+1}^{\pm 1})]= [\rho(b)]$ for all $b\in B_{2g+1}$ ($B_{2g+2}$, resp.)

\bigskip
We shall  pass to a detailed proof  of theorem \ref{thm1} by splitting the argument  into a series of lemmas.
\begin{lem}\label{lem1}
The map $\{\sigma_i\mapsto e_i+1 ~|~ 1\le i\le 2g\}$ 
 ($\{\sigma_i\mapsto e_i+1 ~|~ 1\le i\le 2g+1\}$, resp.)
  defines a representation $\rho: B_{2g+1}\to {\Bbb A}(\mathbf{x}, S_{g,1})$
  ($\rho: B_{2g+2}\to {\Bbb A}(\mathbf{x}, S_{g,2})$, resp.) 
  of the braid group with   an odd  (an even,  resp.)  number of strings 
into a cluster $C^*$-algebra of a Riemann surface with one  cusps
(two cusps, resp.)   
\end{lem}
{\it Proof.}  
We shall prove the case $\rho: B_{2g+1}\to {\Bbb A}(\mathbf{x}, S_{g,1})$ of the braid 
groups with an odd number of strings;  the case of an even number of strings 
is treated likewise.

Let $\{\gamma_1,\dots,\gamma_{2g}\}$ be a chain of simple closed curves
on the surface $S_{g,1}$.  The Dehn twists $\{D_{\gamma_1},\dots,D_{\gamma_{2g}}\}$
around $\gamma_i$ satisfy the braid relations (\ref{eq3}).  The subgroup $SMod~S_{g,1}$
of $Mod~S_{g,1}$ consisting of the automorphisms of $S_{g,1}$ commuting with 
the hyperelliptic involution $\iota$ is isomorphic to the braid group $B_{2g+1}$
(Birman-Hilden's  Theorem \ref{thm2}).  

On the other hand,  each $D_{\gamma_i}$ is a homeomorphism of the Teichm\"uller 
space $T_{g,1}$;  the $D_{\gamma_i}$ induces an inner automorphism of the 
cluster $C^*$-algebra  ${\Bbb A}(\mathbf{x}, S_{g,1})$, see remark \ref{rmk2}. 
We shall denote such an automorphism by
\begin{equation}
\{u_i x u_i^{-1} ~|~ u_i\in {\Bbb A}(\mathbf{x}, S_{g,1}),
~\forall x\in {\Bbb A}(\mathbf{x}, S_{g,1})\}. 
\end{equation}
The units   $\{u_i\in {\Bbb A}(\mathbf{x}, S_{g,1}) ~|~ 1\le i\le 2g\}$ 
satisfy the braid relations:
\begin{equation}\label{eq11}
\cases{
u_iu_{i+1}u_i=u_{i+1}u_i u_{i+1}, &\cr
u_iu_j=u_ju_i,  & if  $|i-j|\ge 2.$}
\end{equation}
Indeed,  from (\ref{eq3}) one gets $SMod~S_{g,1}=\{D_{\gamma_1},\dots, D_{\gamma_{2g}}
~|~ D_{\gamma_i} D_{\gamma_{i+1}}  D_{\gamma_i}=D_{\gamma_{i+1}} D_{\gamma_i} D_{\gamma_{i+1}},
~D_{\gamma_i}D_{\gamma_j}=D_{\gamma_j}D_{\gamma_i} ~\hbox{if} ~|i-j|\ge 2\}$.  
If $Inn~{\Bbb A}(\mathbf{x}, S_{g,1})$ is a group of the inner automorphisms
of the algebra ${\Bbb A}(\mathbf{x}, S_{g,1})$,  then $\{u_i ~|~ 1\le i\le 2g\}$ are generators
of  the  $Inn~{\Bbb A}(\mathbf{x}, S_{g,1})$.  Using the commutative diagram in Figure 1,
one gets from $D_{\gamma_i} D_{\gamma_{i+1}}  D_{\gamma_i}=D_{\gamma_{i+1}} D_{\gamma_i} D_{\gamma_{i+1}}$
the equality   $u_iu_{i+1}u_i=u_{i+1}u_i u_{i+1}$;  similarly, the $D_{\gamma_i}D_{\gamma_j}=D_{\gamma_j}D_{\gamma_i}$
implies the equality $u_iu_j=u_ju_i$  for   $|i-j|\ge 2$.  
In particular, $SMod~S_{g,1}\cong  Inn~{\Bbb A}(\mathbf{x}, S_{g,1})$, where the isomorphism
is given by the formula $D_{\gamma_i}\mapsto u_i$.

It remains to express the units $u_i$ in terms of generators $e_i$ of the algebra  ${\Bbb A}(\mathbf{x}, S_{g,1})$.
The $e_i$ itself is not invertible, but a polynomial $u_i(e_i)=e_i+1$ has an inverse
$u_i^{-1}=-{1\over 2}e_i+1$.  From an isomorphism $B_{2g+1}\cong SMod~S_{g,1}$ 
given by the formula  $\sigma_i\mapsto D_{\gamma_i}$ one gets a representation
$\rho:  B_{2g+1}\to {\Bbb A}(\mathbf{x}, S_{g,1})$ given by the formula 
$\sigma_i\mapsto e_i+1$.  Lemma \ref{lem1} is proved.

\begin{cor}\label{cor1}
The norm closure of a self-adjoint  representation  of  a Temperley-Lieb algebra 
$TL_{2g}({i\sqrt{2}\over 2})$  ($TL_{2g+1}({i\sqrt{2}\over 2})$,  resp.)  is 
a finite-dimensionsal sub-$C^*$-algebra ${\Bbb A}_{2g}$ (${\Bbb A}_{2g+1}$, resp.)   
of  the ${\Bbb A}(\mathbf{x}, S_{g,1})$    (${\Bbb A}(\mathbf{x}, S_{g,2})$, resp.) 
\end{cor}
{\it Proof.}  Let us substitute $u_i=e_i+1$ into the braid relations (\ref{eq11}).  The reader 
is encouraged to verify that relations (\ref{eq11}) are  equivalent to the following system 
of relations:
\begin{equation}\label{eq12}
\cases{
e_i^2=e_i,  ~e_i^*=e_i, &\cr 
e_ie_{i\pm 1}e_i=-2e_i, &\cr
e_ie_j=e_je_i,  & if  $|i-j|\ge 2.$}
\end{equation}
A normalization $e_i'= {i\sqrt{2}\over 2} e_i$ brings (\ref{eq12}) to the form:
\begin{equation}\label{eq13}
\cases{
e_i^2={i\sqrt{2}\over 2} e_i,  &\cr 
e_ie_{i\pm 1}e_i=e_i, &\cr
e_ie_j=e_je_i,  & if  $|i-j|\ge 2.$}
\end{equation}
The relations (\ref{eq13}) are defining relations for a Temperley-Lieb algebra 
$TL_{2g}({i\sqrt{2}\over 2})$  ($TL_{2g+1}({i\sqrt{2}\over 2})$,  resp.),
see e.g.    [Jones 1991,  p. 85]   \cite{J1};   such an algebra is always 
finite-dimensional,  see next lemma. 
Corollary \ref{cor1} follows.

\begin{lem}\label{lem2}
{\bf ([Jones 1991, Section 3.5]   \cite{J1})}
The set ${\cal E} := $ \linebreak
$ \{(e_{i_1}e_{i_1-1}\dots e_{j_1}) \dots  (e_{i_p}e_{i_p-1}\dots e_{j_p}) ~|~ 
1  \le  i_1<\dots <i_p<2g  ~(2g+1, \hbox{resp.}); $
\linebreak
$~1 \le  j_1<\dots <j_p<2g ~(2g+1, \hbox{resp.}); 
~j_1\le i_1, \dots, j_p\le i_p\}$
is multiplicatively closed and 
\begin{equation}
|{\cal E}|\le {1\over n+1}\left(\matrix{2n\cr n}\right),
\end{equation}
where $n=2g$ ($n=2g+1$, resp.)
\end{lem}
{\it Proof.}  An elegant proof of this fact is based on a representation of 
the relations (\ref{eq13}) by the diagrams of the non-crossing strings reminiscent 
of the braid diagrams.

\begin{cor}\label{cor2}
Each   element $e\in {\cal E}$ is equivalent to a projection in the cluster 
$C^*$-algebra ${\Bbb A}(\mathbf{x}, S_{g,1})$   (${\Bbb A}(\mathbf{x}, S_{g,2})$, resp.) 
\end{cor}
{\it Proof.}  
Indeed,  if $e\in {\cal E}$ then the $e^2$ must coincide with one of the elements 
of ${\cal E}$.  But $e^2$ cannot be any such, except for the $e$ itself.  
Thus $e^2=e$, i.e. $e$ is an idempotent.   
On the other hand,  it is well known that each idempotent in a 
$C^*$-algebra is (Murray-von Neumann) equivalent to a projection
in the same algebra,  see e.g.  [Blackadar 1986,  Proposition 4.6.2]  \cite{B}.
Corollary \ref{cor2} follows.

\begin{lem}\label{lem3}
If  $b\in B_{2g+1}$  ($b\in B_{2g+2}$, resp.) is a braid,   there exists 
a Laurent polynomial $[\rho(b)]$  with the integer coefficients  depending on  $2g$ ($2g+1$, resp.) 
variables,   such that $[\rho(b)]\in   K_0({\Bbb A}(\mathbf{x}, S_{g,1}))$   
($[\rho(b)] \in K_0({\Bbb A}(\mathbf{x}, S_{g,2}))$, resp.)
\end{lem}
{\it Proof.}  
We shall prove this fact for  the braid groups with an odd number of strings;  
the case of an even number of strings  is treated likewise.

Let $\{b=\sigma_1^{k_1}\dots \sigma_{2g}^{k_{2g}}\in B_{2g+1} ~|~ k_i\in\mathbf{Z}\}$
be a braid.  By lemma \ref{lem1} such a braid has a representation $\rho(b)$ 
in the cluster $C^*$-algebra  ${\Bbb A}(\mathbf{x}, S_{g,1})$ given by the formula:
\begin{equation}\label{eq15}
\rho(b)=(e_1+1)^{k_1}\dots (e_{2g}+1)^{k_{2g}}\in {\Bbb A}(\mathbf{x}, S_{g,1}).  
\end{equation}
One can unfold the product (\ref{eq15}) into a sum of the monomials in variables 
$e_i$;   by lemma  \ref{lem2}  any such monomial is an element of the set ${\cal E}$. 
In other words, one gets  a finite sum: 
\begin{equation}\label{eq16}
\rho(b)= \left\{\sum_{i=1}^{|{\cal E}|} a_i\varepsilon_i  ~|~ \varepsilon_i\in {\cal E}, ~a_i\in \mathbf{Z}\right\}. 
\end{equation}
(Note that whenever $k_i<0$ the coefficient $a_i$ is a rational number, 
but clearing the denominators we can assume  $a_i\in \mathbf{Z}$.)  
On the other hand,  corollary \ref{cor2}  says that each $\varepsilon_i$ is a
projection;  therefore $\varepsilon_i$ defines an equivalence class 
$[\varepsilon_i]\in   K_0({\Bbb A}(\mathbf{x}, S_{g,1}))$ of projections in 
the  cluster $C^*$-algebra  ${\Bbb A}(\mathbf{x}, S_{g,1})$.  Thus one gets 
$[\rho(b)] \in  K_0({\Bbb A}(\mathbf{x}, S_{g,1}))$ given by a finite sum
\begin{equation}\label{eq17}
[\rho(b)]= \left\{\sum_{i=1}^{|{\cal E}|} a_i [\varepsilon_i]  ~|~ [\varepsilon_i]\in K_0({\Bbb A}(\mathbf{x}, S_{g,1})), 
~a_i\in \mathbf{Z}\right\}. 
\end{equation}
But  $K_0({\Bbb A}(\mathbf{x}, S_{g,1}))\cong {\cal  A}(\mathbf{x}, S_{g,1})\subset 
\mathbf{Z}[\mathbf{x}^{\pm 1}]$;  in particular,  $[\rho(b)]\in \mathbf{Z}[\mathbf{x}^{\pm 1}]$
is a Laurent polynomial with the integer coefficients.  

To calculate the number of variables in  $[\rho(b)]$,  recall that 
$rank~ {\cal  A}(\mathbf{x}, S_{g,1})=6g-3$;  on the other hand,
a fundamental domain of the Riemann surface $S_{g,1}$ is a paired 
$(4g+2)$-gon,  where one pair of sides  corresponds to a boundary
component obtained from the cusp, see Section 2.1.
Since the boundary component contracts to a cusp,  one gets 
a paired $4g$-gon whose triangulation requires $4g-3$ interior geodesic arcs. 
Thus the cluster $|\mathbf{x}|=6g-3$ can be written in the form:
\begin{equation}\label{eq18}
\mathbf{x}=(x_1,\dots,x_{2g}; y_1,\dots, y_{4g-3}),
\end{equation}
where $x_i$ are  mutable and $y_i$ are  frozen  
variables [Williams 2014,  Definition 2.6] \cite{Wil1}.   
One can always assume $y_i=Const$  and therefore 
the Laurent polynomial $[\rho(b)]$ depends on the $2g$ variables $x_i$.
Lemma \ref{lem3} follows.

\begin{lem}\label{lem4}
The Laurent polynomial $[\rho(b)]$ is a topological invariant of the closure of 
the  braid $b\in B_{2g+1}$ ($b\in B_{2g+2}$, resp.)  
\end{lem}
{\it Proof.}  
Again we shall prove the case $b\in B_{2g+1}$;  the case  $b\in B_{2g+1}$
can be treated similarly and is left to the reader. To prove that  $[\rho(b)]$
is a topological invariant, it is enough to demonstrate that:

\medskip
(i)   $[\rho(gbg^{-1})]= [\rho(b)]$ for all $g\in B_{2g+1}$;

\smallskip
(ii)  $[\rho(b\sigma_{2g+1}^{\pm 1})]= [\rho(b)]$ for the generator  $\sigma_{2g+1}\in B_{2g+2}$.

\medskip
(i)  Recall that the $K$-theory of an $AF$-algebra ${\Bbb A}(\mathbf{x}, S_{g,1})$ 
can be recovered from the canonical trace $\tau:  {\Bbb A}(\mathbf{x}, S_{g,1})\to \mathbf{C}$
 [Blackadar 1986,  Section 7.3]  \cite{B}. On the other hand, the trace $\tau$ is a
 character of the representation $\rho: B_{2g+1}\to  {\Bbb A}(\mathbf{x}, S_{g,1})$; 
 in particular, $\tau(\rho(g)\rho(b)\rho(g^{-1}))=\tau(\rho(b))$. 
 In other words, $[\rho(gbg^{-1})]= [\rho(b)]$ for all $g\in B_{2g+1}$.
 Item (i) follows.  
 
 \medskip
 (ii)  Let ${\Bbb A}_{2g}$ be a finite-dimensional $C^*$-algebra of corollary \ref{cor1}.
 Let   $u_{2g+1}$ and $u_{2g+1}^{-1}$   be a unit and its inverse 
in the algebra ${\Bbb A}(\mathbf{x}, S_{g,1})$  given by the formulas:
\begin{equation}\label{eq19}
u_{2g+1}=u_{2g+1}^{-1}:=2e_{2g+1}-1.  
\end{equation}
Clearly the units $u_{2g+1}^{\pm 1}\not\in {\Bbb A}_{2g}$ and they make the
group $G\cong \mathbf{Z}/2\mathbf{Z}$ under a multiplication. 
We shall consider an extension ${\Bbb A}_{2g}\rtimes G$ of the algebra 
${\Bbb A}_{2g}$ given by the formal sums
\begin{equation}\label{eq20}
{\Bbb A}_{2g}\rtimes G:= \left\{\sum_{\gamma\in G} a_{\gamma} u_{\gamma} 
~|~ u_{\gamma}\in {\Bbb A}_{2g}\right\}. 
\end{equation}
(The algebra of formal  sums (\ref{eq20}) is isomorphic to a crossed product 
$C^*$-algebra of the algebra ${\Bbb A}_{2g}$ by the outer automorphisms 
$\alpha$ given the elements $u_{\gamma}\in \{u_{2g+1}^{\pm 1}\}$;  hence 
our notation.)  It is well known that 
\begin{equation}\label{eq21}
{\Bbb A}_{2g}\rtimes G\cong M_2({\Bbb A}_{2g}^{\alpha}),
 \end{equation}
 where ${\Bbb A}_{2g}^{\alpha}\subset {\Bbb A}_{2g}$ is a fixed-point 
 algebra of the  automorphism $\alpha: {\Bbb A}_{2g}\to {\Bbb A}_{2g}$,
 see e.g. [Fillmore 1996, Section 3.8.5]  \cite{F}.  
By the second of  formulas (\ref{eq11}),  one gets $u_i u_{2g+1}=
u_{2g+1}u_i$ for all $1\le i\le 2g-1$;   in other words,
\begin{equation}\label{eq22}
u_i=u_{2g+1}u_i u_{2g+1}^{-1}
 \end{equation}
for all $1\le i\le 2g-1$.  Since all generators $u_i$ of algebra ${\Bbb A}_{2g-1}$ 
 are fixed by the automorphism $\alpha$,  one gets an isomorphism 
\begin{equation}\label{eq23}
{\Bbb A}_{2g}^{\alpha}\cong {\Bbb A}_{2g-1}. 
 \end{equation}
On the other hand, $K_0(M_2(A))\cong K_0(A)$ by a stability of the 
$K$-theory  [Blackadar 1986, Section 5.1]  \cite{B};  thus formulas 
(\ref{eq21}) and (\ref{eq23}) imply an isomorphism:
\begin{equation}\label{eq24}
 K_0({\Bbb A}_{2g}\rtimes G)\cong K_0({\Bbb A}_{2g-1}). 
 \end{equation}
It remains to notice that if one maps the generators $\sigma_{2g+1}^{\pm 1}
\in B_{2g+2}$  into the units $u_{2g+1}^{\pm 1}\in {\Bbb A}_{2g}\rtimes G$,
then formulas (\ref{eq20}) and (\ref{eq24}) imply the equality
\begin{equation}\label{eq25}
 [\rho(b\sigma_{2g+1}^{\pm 1})]= [\rho(b)]
 \end{equation}
for all $b\in B_{2g}$.  The item (ii) follows from (\ref{eq25})  and lemma 
\ref{lem4} is proved.

\bigskip
Theorem \ref{thm1} follows from lemmas \ref{lem1},  \ref{lem3} and \ref{lem4}.

\section{Examples}
To illustrate theorem \ref{thm1}, we shall consider a representation 
$\rho: B_2\to {\Bbb A}(\mathbf{x}, S_{0,2})$  ($\rho: B_3\to {\Bbb A}(\mathbf{x}, S_{1,1})$, resp.);
it will be shown that for such a representation the Laurent polynomials $[\rho(b)]$ correspond 
to the Jones (HOMFLY, resp.)  invariants of knots and links.

\subsection{Jones polynomials}
If $g=0$ and $n=2$, then $S_{0,2}$ is a sphere with two cusps; 
the  $S_{0,2}$  is homotopy equivalent to an annulus
${\goth A}:=\{z=u+iv\in\mathbf{C} ~|~ r\le |z|\le R\}$.  
The Riemann surface ${\goth A}$ has an ideal triangulation $T$  with 
one marked point on each boundary component  given by the matrix:
\begin{equation}\label{eq26}
B_T=\left(\matrix{0 & 2\cr -2 & 0}\right),   
\end{equation}
see [Fomin,  Shapiro  \& Thurston  2008, Example 4.4]  \cite{FoShaThu1}. 
Using the exchange relations (\ref{eq1})  the reader  can verify that the 
cluster $C^*$-algebra ${\Bbb A}(\mathbf{x}, S_{0,2})$ is given by the 
Bratteli diagram shown in Figure 2;   the ${\Bbb A}(\mathbf{x}, S_{0,2})$
coincides with the so-called {\it GICAR algebra}  [Bratteli 1972, Section 5.5]  \cite{Bra1}.

\begin{figure}[here]
\begin{picture}(100,100)(-150,150)

\put(50,200){\circle{3}}

\put(33,188){\circle{3}}
\put(67,188){\circle{3}}

\put(50,177){\circle{3}}
\put(16,177){\circle{3}}
\put(84,177){\circle{3}}

\put(-1,164){\circle{3}}
\put(34,164){\circle{3}}
\put(68,164){\circle{3}}
\put(103,164){\circle{3}}


\put(49,199){\vector(-3,-2){15}}
\put(51,199){\vector(3,-2){15}}

\put(32,187){\vector(-3,-2){15}}
\put(34,187){\vector(3,-2){15}}

\put(66,187){\vector(-3,-2){15}}
\put(68,187){\vector(3,-2){15}}


\put(14,175){\vector(-3,-2){15}}
\put(17,175){\vector(3,-2){15}}

\put(49,175){\vector(-3,-2){15}}
\put(51,175){\vector(3,-2){15}}

\put(83,175){\vector(-3,-2){15}}
\put(86,175){\vector(3,-2){15}}

\put(-10,155){$\dots$}
\put(27,155){$\dots$}
\put(64,155){$\dots$}
\put(101,155){$\dots$}

\end{picture}
\caption{Bratteli diagram of the algebra ${\Bbb A}(\mathbf{x}, S_{0,2})$.} 
\end{figure}
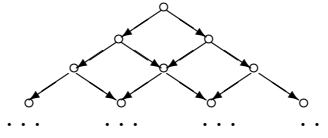

\bigskip\noindent
The cluster $\mathbf{x}=(x; c)$ consists of a mutable variable $x$ and 
a coefficient  $c\in ({\Bbb P}, \oplus, \bullet)$.  Theorem \ref{thm1} 
says that there exists a representation   
\begin{equation}\label{eq27}
 \rho: B_2\to {\Bbb A}(\mathbf{x}, S_{0,2}),
 \end{equation}
such that $[\rho(b)]\in\mathbf{Z}[x^{\pm 1}]$ is a topological invariant 
of the closure $L$ of $b\in B_2$;  since  the Laurent  polynomial $[\rho(b)]$
depends on $x$ and  $c$,  we shall write it  $[\rho(b)](x,c)$.   Let $N\ge 1$ be the minimal 
number  of the overpass (underpass, resp.) crossings added 
to the unknot $K$ to get the link $L$  [Jones 1985, Figure 2]  \cite{Jon1}.  
The following result compares
the $[\rho(b)](x,c)$  with the Jones polynomial $V_L(t)$.  
\begin{cor}\label{cor3}
\begin{equation}\label{jones}
V_L(t)=\left(-{\sqrt{t} \over t+1}\right)^N ~[\rho(b)](t, -t^2).
\end{equation}  
\end{cor}
{\it Proof.} 
Recall that each  polynomial $[\rho(b)](x,c)$ is  obtained from an initial seed $(\mathbf{x}, B_T)$
by a finite number of mutation given by the exchange relations (\ref{eq1});
likewise,   each polynomial $V_L(t)$ can be obtained from the $V_K(t)=1$
 using  the skein relation (\ref{skein}).
Roughly speaking, the idea is to show that (\ref{eq1}) and    (\ref{skein})  are equivalent relations
up to  a   multiple    $-{\sqrt{t}\over t+1}$.
Indeed,  consider the Laurent polynomials: 
\begin{equation}\label{eq28} 
\cases{ W_{L^+}=\left(-{t+1\over\sqrt{t}}\right) ~V_{L^+}&\cr
W_{L^-}=\left(-{t+1\over\sqrt{t}}\right) ~V_{L^-} .&}
\end{equation}
The skein relation  (\ref{skein})  for the $W_{L^{\pm}}$  takes the form:
\begin{equation}\label{eq29} 
V_L={t^2\over t^2-1} W_{L^+}-{1\over t^2-1} W_{L^-}. 
\end{equation}
The substitution 
\begin{equation}\label{eq30} 
\cases{ 
V_L=x_k'&\cr
W_{L^+}=  {1\over x_k}  \prod_{i=1}^n  x_i^{\max(b_{ik}, 0)}&\cr
W_{L^-}= {1\over x_k}  \prod_{i=1}^n  x_i^{\max(-b_{ik}, 0)}  &\cr
c_k=-t^2&}
\end{equation}
transforms the skein relation (\ref{eq29}) to the exchange relations
(\ref{eq1}).   It remains to observe from (\ref{eq28}),  that an extra crossing added to the
unknot $K$ corresponds to  a multiplication of   $[\rho(b)](x,c)$ by $-{\sqrt{t}\over t+1}$;
the minimal number $N$ of such crossings required to get the link  $L$ from $K$ gives
us   the $N$-th power of  the multiple.
Corollary \ref{cor3} follows.

\subsection{HOMFLY polynomials}
If $g=n=1$,  then $S_{1,1}$ is a torus with a cusp. The matrix $B_T$ associated 
to an ideal triangulation of the Riemann surface $S_{1,1}$    has the form: 
 \begin{equation}\label{eq31}
 B_T=\left(
 \matrix{0 & 2 & -2\cr
              -2 & 0 & 2\cr
               2 & -2 & 0}
              \right),
 \end{equation}
see [Fomin,  Shapiro  \& Thurston  2008, Example 4.6]  \cite{FoShaThu1}. 
Using the exchange relations (\ref{eq1})  the reader  can verify that the 
cluster $C^*$-algebra ${\Bbb A}(\mathbf{x}, S_{1,1})$ is given by the 
Bratteli diagram shown in Figure 3;   the ${\Bbb A}(\mathbf{x}, S_{1,1})$
coincides with the  {\it Mundici  algebra} ${\goth M}_1$  [Mundici 1988]  \cite{Mun1}. 

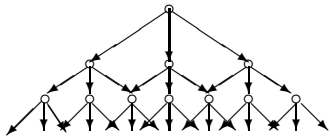
\begin{figure}[here]
\begin{picture}(100,100)(-160,130)

\put(50,200){\circle{3}}
\put(50,179){\circle{3}}
\put(20,179){\circle{3}}
\put(80,179){\circle{3}}


\put(3,166){\circle{3}}
\put(20,166){\circle{3}}
\put(36,166){\circle{3}}
\put(50,166){\circle{3}}
\put(65,166){\circle{3}}

\put(80,166){\circle{3}}
\put(98,166){\circle{3}}


\put(49,199){\vector(-3,-2){29}}
\put(51,199){\vector(3,-2){29}}
\put(50,200){\vector(0,-1){20}}


\put(19,178){\vector(-3,-2){15}}
\put(21,178){\vector(3,-2){15}}
\put(20,178){\vector(0,-1){10}}


\put(50,178){\vector(-3,-2){15}}
\put(50,178){\vector(3,-2){15}}
\put(50,178){\vector(0,-1){10}}


\put(79,178){\vector(-3,-2){15}}
\put(81,178){\vector(3,-2){15}}
\put(80,178){\vector(0,-1){10}}


\put(1.5,165){\vector(-1,-1){13}}
\put(4,164){\vector(2,-3){7}}
\put(2.5,164){\vector(0,-1){10}}

\put(19,164){\vector(-1,-1){10}}
\put(21,164){\vector(1,-1){10}}
\put(20,164){\vector(0,-1){10}}

\put(49,164){\vector(-1,-1){10}}
\put(51,164){\vector(1,-1){10}}
\put(50,164){\vector(0,-1){10}}

\put(79,164){\vector(-1,-1){10}}
\put(81,164){\vector(1,-1){10}}
\put(80,164){\vector(0,-1){10}}

\put(98,164){\vector(-1,-1){10}}
\put(100,164){\vector(1,-1){10}}
\put(98,164){\vector(0,-1){10}}


\put(65,164){\vector(-1,-1){10}}
\put(65,164){\vector(1,-1){10}}
\put(65,164){\vector(0,-1){10}}

\put(36,164){\vector(-1,-1){10}}
\put(36,164){\vector(1,-1){10}}
\put(36,164){\vector(0,-1){10}}


\end{picture}
\caption{Bratteli diagram of the algebra ${\Bbb A}(\mathbf{x}, S_{1,1})$.} 
\end{figure}

\bigskip\noindent
The formula (\ref{eq18}) implies  that  cluster $\mathbf{x}=(x_1,x_2;  y_1)$ consists of two mutable 
variables $x_1, x_2$ and a frozen variable $y_1$.    Theorem \ref{thm1} 
says that there exists a representation:   
\begin{equation}\label{eq32}
 \rho: B_3\to {\Bbb A}(\mathbf{x}, S_{1,1}),
 \end{equation}
such that $[\rho(b)]\in\mathbf{Z}[x_1^{\pm 1}, x_2^{\pm 1}]$ is a topological invariant 
of the closure $L$ of $b\in B_3$;  since  the Laurent  polynomial $[\rho(b)]$
depends on two variables $x_1, x_2$ and  two coefficients  $c_1,c_2\in  ({\Bbb P}, \oplus, \bullet)$,  
we shall write it  $[\rho(b)](x_1, x_2;  c_1, c_2)$. 
The following result says that the  $[\rho(b)](x_1, x_2; c_1, c_2)$  for  special values  of $c_i$ 
is related to the HOMFLY  polynomial $\rho_L(l,m)$.  (Unlike (\ref{jones}),  an explicit formula
for such a relationship seems to be complicated.) 
\begin{cor}\label{cor4}
The exchange relations (\ref{eq1})  with  matrix  $B$ given by  (\ref{eq31})
imply the skein relation (\ref{homfly}) for the HOMFLY polynomial   $\rho_L(l,m)$.
\end{cor}
{\it Proof.}    Since the variable $y_1$ is frozen,  we  consider
a reduced matrix: 
\begin{equation}\label{eq33}
\tilde B_T=\left(\matrix{0 & 2\cr -2 & 0}\right)
\end{equation}
and  our seed has the form $(\mathbf{x},  \tilde B_T)$,
where $\mathbf{x}=(x_1, x_2;  c_1, c_2)$.  
The exchange relations (\ref{eq1}) for the variables $x_3, x_4, x_5$
and the coefficient $c_3$  imply  the following system of equations: 
\begin{equation}\label{eq34} 
\cases{ 
x_3={c_1+x_2^2\over (c_1+1)x_1}&\cr
x_4={c_2 x_3^2+1\over (c_2+1)x_2}&\cr
x_5={c_3+x_4^2\over (c_3+1)x_3}&\cr
c_3={1\over c_1}.&}
\end{equation}
Clearing the denominators in (\ref{eq34}), one gets an expression:
\begin{equation}\label{eq35}
c_1x_3+c_3x_5={c_1+x_2^2-x_1x_3\over x_1} + {c_3+x_4^2-x_3x_5\over x_3}. 
\end{equation}
We exclude $x_4= {c_2 x_3^2+1\over (c_2+1)x_2}$  and $c_3={1\over c_1}$ 
in (\ref{eq35}) and get  the following   equality of the Laurent polynomials: 
\begin{equation}\label{eq36}
c_1x_3+{1\over c_1} x_5+c_2W=0,
\end{equation}
where $W:= -c_2\left[  {c_1+x_2^2-x_1x_3\over c_2^2 x_1}+
{c_1^{-1}-x_3x_5\over c_2^2x_3}+
{1\over x_3}\left({c_2x_3^2+1\over c_2(c_2+1)x_2}\right)^2\right]$. 
The substitution:
\begin{equation}\label{eq37} 
\cases{ 
c_1=x_1=l&\cr
c_2= x_2=m&\cr
x_3= \rho_{L^+} &\cr
x_5=\rho_{L^-} &\cr
W=\rho_L&}
\end{equation}
brings (\ref{eq36})  to the skein relation (\ref{homfly}). 
Corollary \ref{cor4} follows.

\bigskip\noindent
{\sf Acknowledgment.} 
I thank the members of the SAG group (Ibrahim Assem, Thomas Br\"ustle,  Virginie Charette, 
Tomasz Kaczynski,  Shiping Liu  and Vasilisa Shramchenko)  at the Department of Mathematics 
of the University of Sherbrooke for their interest, hospitality and excellent 
working conditions.



\vskip1cm

\textsc{The Fields Institute for Research in Mathematical Sciences, Toronto, ON, Canada,  
E-mail:} {\sf igor.v.nikolaev@gmail.com}


\end{document}